\documentclass{amsart}
\usepackage{amsmath,amsfonts,amsthm,amssymb}
\usepackage[alphabetic]{amsrefs}

\newtheorem{lem}{Lemma}[section]
\newtheorem{de}[lem]{Definition}
\newtheorem{tw}[lem]{Theorem}

\newtheorem{cor}[lem]{Corollary}
\newtheorem{prop}[lem]{Proposition}

\newcommand {\mr}{\mathrm}
\newcommand {\lk}{\left\{}
\newcommand {\rk}{\right\}}

\newcommand {\ov}{\overline}
\newtheorem*{gtw}{Main Theorem}
\newcommand {\rems}{\noindent {\bf Remarks. }}
\newcommand {\rem}{\noindent {\bf Remark. }}

\usepackage[dvips]{graphicx}

\begin {document}

\title{A construction of hyperbolic Coxeter groups}
\author{Damian Osajda}
\date{\today}
\thanks{Partially supported by MNiSW
grant N201 012 32/0718 and by the ANR grant "Cannon". This research was supported through the programme "Oberwolfach Leibniz Fellows" by the Mathematisches Forschungsinstitut Oberwolfach in 2010}
\address{Instytut Matematyczny,
Uniwersytet Wroc\l awski\\
pl. Grunwaldzki 2/4,
50--384 Wroc{\l}aw, Poland}
\email{dosaj@math.uni.wroc.pl}

\begin{abstract}
We give a simple construction of Gromov hyperbolic Coxeter groups of arbitrarily large virtual cohomological dimension. Our construction provides new examples of such groups. Using this one can construct e.g. new groups having some interesting asphericity properties.
\end{abstract}

\maketitle

\section{Introduction}
\label{intro}
The question of constructions of highly dimensional Gromov hyperbolic groups was raised several times in the past.
Gromov's loose conjecture \cite{G} (see its precise version on the Bestvina's problem list \cite{Be-p}) stated that
a construction of such groups always involves nontrivial number theoretic tools (see also the discussion in \cite{JS1}).
Eventually Januszkiewicz-\'Swi\c atkowski \cite{JS1} provided a geometric
construction of Gromov hyperbolic Coxeter groups in every dimension.
Later several similar constructions appeared \cites{H,JS2,AB+}.
All of them are using a fairly advanced machinery of complexes of groups and, moreover, the groups obtained those ways are always systolic; cf. \cite{JS2}. As shown by Januszkiewicz-\'Swi\c atkowski \cite{JS3} and the author \cites{O-ciscg,O-ib7scg} systolic groups satisfy some very restrictive asphericity properties that make them in a way asymptotically two-dimensional. In particular they do not "contain asymptotically" spheres of dimension two and more.

In this paper we give a simple geometric construction of Gromov hyperbolic Coxeter groups of arbitrarily large virtual cohomological dimension.
Actually, it is the simplest construction of highly dimensional Gromov hyperbolic groups known to us.
Our method is quite elementary and uses only (in the simplest version) basic facts about right-angled Coxeter groups. It allows us to construct highly dimensional Gromov hyperbolic Coxeter groups that are not systolic---they may contain spheres at infinity. We are also able to provide new examples of systolic groups, highly dimensional asymptotically hereditarily aspherical (shortly AHA) groups (cf. \cite{JS3}) and groups
with strongly hereditarily aspherical (SHA) Gromov boundaries (cf. \cite{O-ib7scg}). This paper bases on tools and ideas around weakly systolic complexes introduced and developed by the author in \cite{O-sdn}.
In particular our construction gives examples of weakly systolic groups of arbitrarily large dimension that are not systolic.

It should be noticed here that for some time it was believed---cf. a conjecture of Moussong \cite{M}---that there is a universal bound on the virtual cohomological dimension of any Gromov hyperbolic Coxeter group.
One of motivations for that  was the result of Vinberg \cite{V} stating that if a reflection group acts cocompactly on an $n$-dimensional real hyperbolic space $\mathbb H^n$ then $n\leq 29$.
Eventually Januszkiewicz-\'Swi\c atkowski \cite{JS1} disproved it, but our examples are the first non-systolic examples of Gromov hyperbolic Coxeter groups in high dimensions.
\medskip

The construction is as follows (cf. the following sections for explanations of all the notions involved). To construct a Gromov hyperbolic right angled Coxeter group all we need is a finite $5$-large (or flag no-square) simplicial complex. We construct such simplicial complexes inductively. Their "dimension" increases each time. The induction step is itself
divided into the following three steps.

\vspace{0.5cm}
{\bf The basic construction.}
\medskip

\noindent
{\bf Step 1.}
Let $X$ be a finite $5$-large (i.e. flag no-square) simplicial complex with $H^n(X)\neq 0$.
Let $(W,S)$ be the Coxeter system whose nerve is $X$, i.e. $X=L(W,S)$.
Then the virtual cohomological dimension of $W$ is $\mr{vcd}\; W\geq n+1$.
The complex $X$ is the link of every vertex of the Davis complex $\Sigma=\Sigma(W,S)$ of $(W,S)$.
\medskip

\noindent
{\bf Step 2.}
Choose a torsion-free subgroup $H$ of sufficiently large finite index in $W$. Then $Y=\Sigma/ H$ is a locally $5$-large cubical complex (i.e. with links being flag no-square simplicial complexes) with $H^{n+1}(Y)\neq 0$.
\medskip

\noindent
{\bf Step 3.}
Let $X'=Th(Y)$ be the simplicial thickening of $Y$ (i.e. the simplicial complex obtained by replacing cubes by simplices spanned on their vertices). Then $X'$ is a finite $5$-large simplicial complex with $H^{n+1}(X')\neq 0$.

\vspace{0.5cm}
This is the end of the inductive step. The constructed complex $X'$ can be now used as the initial complex $X$ again in Step 1. The crucial fact for this construction to work and the main result of the paper is the following.

\begin{gtw}[cf. Theorem \ref{maint} in Section \ref{proof}]
\label{main}
Let $X$ be a finite $5$-large simplicial complex such that $H^n(X)\neq 0$.
Then there exists a torsion-free finite index subgroup $H$ of $W$ such that the complex $X'$ obtained by the basic construction above is a finite $5$-large simplicial complex with $H^{n+1}(X')\neq 0$. Thus the right-angled Coxeter group with the nerve $X'$ is Gromov hyperbolic of virtual cohomological dimension at least $n+2$.
\end{gtw}

Starting with a finite $5$-large $n_0$-dimensional simplicial complex $X_0$, after performing $k$ times Steps 1--3, we get a finite $5$-large simplicial complex $X$ with $H^{n_0+k}(X)\neq 0$. Thus the right-angled Coxeter group with the nerve $X$ is Gromov hyperbolic of virtual cohomological dimension at least $(n_0+k+1)$.
\medskip

The main idea behind this construction is that one uses the finite quotient $X$ of an "$n$-dimensional" complex (of the Davis complex) as the link of vertices in a new complex (the new Davis complex). Thus the "dimension" of the new complex jumps up at least one (since the new complex is a union of cones over its links). Then one proceeds inductively.

After the description of the construction given above an expert reader can actually skip Sections \ref{pre}--\ref{final}. In the remaining part of the introduction we describe the content of next sections pointing out other results related to the main version of the construction.

Sections \ref{pre} and \ref{prel-sd} are preliminary sections. In particular, in Section \ref{prel-sd} we reprove some technical results from \cite{O-sdn} concerning thickenings of cubical complexes.

In Section \ref{proof} we present in details the basic construction following the scheme given above. In particular we prove in details Main Theorem above; cf. Theorem \ref{maint}.

Results from Section \ref{new} show that our construction provides examples of Gromov hyperbolic groups of arbitrarily large dimension that are asymptotically different from the groups known before. In particular we show that our groups can "contain asymptotically" spheres of dimension up to $3$; cf. Corollary \ref{nonsys} and Examples, Remarks afterwards.

On the other hand, in Section \ref{newa} we show that we can also produce systolic groups (cf. Proposition \ref{syst}) and weakly systolic groups (not systolic a priori), which asymptotically look very similar to systolic ones (cf. Proposition \ref{sd*}, Corollary \ref{asph} and Remarks afterwards).

Finally, in the concluding Section \ref{final} we give some final remarks, e.g. about variants of the basic construction.
\medskip

\noindent
{\bf Acknowledment.} I thank Jan Dymara, Fr\'ed\'eric Haglund, Tadeusz Januszkiewicz and Jacek \'Swi\c atkowski for helpful conversations.

During writing of this paper I found out that the same idea of the basic construction was discovered independently and at the same time by Fr\'ed\'eric Haglund (unpublished).

\section{Preliminaries}
\label{pre}
\subsection{Simplicial complexes}
\label{simpl}

Let $X$ be a simplicial complex. The $i$-skeleton of $X$ is denoted by $X^{(i)}$.
A subcomplex $Y$ of $X$ is \emph{full} if
every subset $A$ of vertices of $Y$ contained in a simplex of $X$, is
contained in a simplex of $Y$.
For a finite set $A=\lk v_1,...,v_k \rk $ of vertices of $X$, by $\mr{span}(A)$ or by $\langle v_1,...,v_k \rangle$ we denote the \emph{span} of $A$, i.e. the smallest full subcomplex of $X$ containing $A$.
A simplicial complex $X$ is \emph{flag} whenever every finite set of vertices of $X$ joined pairwise by edges in $X$, is contained in a simplex of $X$.
A \emph{link} of a simplex $\sigma$ of $X$ is a simplicial complex $X_{\sigma}=\lk \tau | \; \tau \in X \; \& \; \tau \cap \sigma=\emptyset \; \& \; \mr{span}(\tau \cup \sigma)\in X \rk$.

Let $k\geq 4$. A \emph{$k$-cycle} $(v_0,...,v_{k-1},v_0)$ is a triangulation of a circle consisting of $k$ edges ($\langle v_i,v_{i+1\; (\mr{mod}\; k)}\rangle $) and $k$ vertices: $v_0,...,v_{k-1}$.
A \emph{$k$-wheel (in $X$)} $(v_0;v_1,...,v_k)$ (where $v_i$'s are vertices of
$X$) is a full subcomplex of $X$ such that $(v_1,...,v_k,v_1)$ is a
$k$-cycle and $v_0$ is joined (by an edge in $X$) with $v_i$, for $i=1,...,k$.
A \emph{$k$-wheel with a pendant triangle (in $X$)} $(v_0;v_1,...,v_k;t)$ is a
subcomplex of $X$ being the union of a $k$-wheel $(v_0;v_1,...,v_k)$ and a triangle
$\langle v_1,v_2,t \rangle \in X$, with $t\neq v_i$, $i=0,1,...,k$.
for $k\geq 4$, a flag simplicial complex $X$ is \emph{$k$-large} if there are no $j$-cycles being full subcomplexes of $X$, for $j<k$ ($4$-large means simply flag). The term \emph{flag no square} is sometimes used instead of $5$-large.
$X$ is \emph{locally $k$-large} if all its links are $k$-large.

\subsection{Cubical complexes}
\label{cc}
\emph{Cubical complexes} are cell complexes in which every cell is isomorphic to a standard (Euclidean or hyperbolic) cube; see \cite[Chapter I.7]{BrHa} and \cite[Appendix A]{Dav-b} for a precise definition. It means in particular that two cubes in a cubical complex intersect along a single subcube. The \emph{link} $Y_v$ of a vertex $v$ of a cubical complex $Y$ is a simplicial complex that can be identified with a small sphere around $v$ (simplices of $Y_v$ are intersections of the sphere with cubes).
A cubical complex is \emph{locally $k$-large} (resp. \emph{locally flag}) if links of its vertices are $k$-large (resp. flag).

A lemma of Gromov (cf. \cite[Appendix I]{Dav-b}) states that a simply connected locally flag (resp. locally $5$-large) cubical complex admits a metric of non-positive (resp. negative) curvature, or CAT(0) (resp. CAT(-1)) metric.

\subsection{Coxeter groups}
\label{cox}
We use terminology and notations from the Davis' book \cite{Dav-b}.
A \emph{Coxeter group} is given by a presentation $W=\langle S| (st)^{m_{st}}; \; s,t\in S\rangle$, where $S$ is a finite set, $m_{st}\in \mathbb{N}\cup \lk \infty \rk$, $m_{st}=m_{ts}$ and $m_{st}=1$ iff $t=s$ (here $(st)^{\infty}$ means no relation). A Coxeter group (or a \emph{Coxeter system} $(W,S)$) is called \emph{right-angled} if $m_{st}\in \lk 1,2,\infty \rk$.
A \emph{special subgroup} $W_T$ of a Coxeter group $W$ is a subgroup generated by a subset $T\subseteq S$. A subset $T\subseteq S$ is called \emph{spherical} if $W_T$ is finite. In that case $W_T$ is called also \emph{spherical}. By $\mathcal S$ we denote the poset (wrt inclusions) of spherical subsets of $S$. Its geometric realization is denoted by $K$. The poset of all nonempty spherical subsets is an abstract simplicial complex: the \emph{nerve} $L=L(W,S)$ of the Coxeter system $(W,S)$. For $T\subseteq \mathcal S$, by $\sigma_T$ we denote the simplex of $L$ spanned by $T$ ($\sigma_T=\emptyset$ if $T=\emptyset$).
The geometric realization of the poset (wrt inclusions) $\bigcup_{T\in \mathcal{S}} W/W_T$ is called the \emph{Davis complex} and is denoted by $\Sigma=\Sigma(W,S)$.
In the right-angled case $\Sigma$ possess a natural structure of a locally flag cubical complex.
For $s\in S$ we define $K_s$ as the union of the simplices in $K$ with minimum vertex $\lk s \rk$. For $T\subseteq S$ we define $K^T=\bigcup_{s\in T} K_s$. It can be shown that $L-\sigma_T$ deformation retracts onto $K^{S-T}$; cf. \cite[Lemma A.5.5]{Dav-b}.

Moussong \cite{M} showed the following.
\begin{tw}(Hyperbolic right-angled Coxeter group)
\label{mous}
A right-angled Coxeter group $(W,S)$ is Gromov hyperbolic iff its nerve $L(W,S)$
is a $5$-large (i.e. flag no-square) simplicial complex.
\end{tw}

In fact in that case the Davis complex $\Sigma(W,S)$ possess a natural structure of a locally $5$-large (i.e. CAT(-1); cf. Section \ref{cc}) cubical complex.

\subsection{Virtual cohomological dimension}
\label{vcd}
Recall (cf. \cite[Chapter 8.5]{Dav-b}) that the \emph{cohomological dimension}
of a group $G$ is defined as
\begin{align*}
\label{cdform}
\mr{cd}\; G=\sup \lk n|\; H^n(G;M)\neq 0 \; \mr{for\; some} \;
\mathbb ZG\mr{-module} \; M \rk.
\end{align*}
If $G$ has nontrivial torsion then $\mr{cd}\; G=\infty$.
Thus for \emph{virtually torsion-free} groups (i.e. groups having a torsion-free subgroup of finite index) the following notion of dimension is more convenient.

The \emph{virtual comhomological dimension} of a group $G$, denoted $\mr{vcd}\; G$, is the cohomological dimension of any torsion-free finite index subgroup of $G$.
For a Coxeter system $(W,S)$ as above we have the following.

\begin{tw}[{\cite[Corollary 8.5.5]{Dav-b}}]
\label{vcdcox}
\begin{align*}
\mr {vcd}\; W & =\max \lk n|\; H^{n}( K,K^{S-T})\neq 0, \; \mr{for \; some}\; T\in \mathcal S \rk\\
& = \max \lk n|\; \ov H^{n-1}( L- \sigma_T)\neq 0, \; \mr{for \; some}\; T\in \mathcal S \rk.
\end{align*}
\end{tw}

\subsection{Weakly systolic complexes}
\label{weak}

Recall---cf. \cite{JS1} and \cite{Ch4} (where the name \emph{bridged complexes} is used)---that a flag simplicial complex $X$ is called \emph{$k$-systolic}, $k\geq 6$, if it is simply connected and locally $k$-large.

In \cite{O-sdn}, the following weakly systolic complexes were introduced as a generalization of
systolic complexes (see also \cite{OCh}).

\begin{de}[Weakly systolic complex]
\label{weaks}
A flag simplicial complex $X$ is \emph {weakly systolic} if for every vertex $v$ of $X$
and for every $i=1,2,...$ the following condition holds.
For every simplex $\sigma \subseteq S_{i+1}(v,X)$ the intersection $X_{\sigma}\cap B_i(v,X)$ is a non-empty simplex.
\end{de}

Systolic complexes are weakly systolic \cite{O-sdn}. By Lemma \ref{5l thick} below, thickenings of locally $5$-large (or CAT(-1); cf. Section \ref{cc}) cubical complexes are weakly systolic.Thus the class of \emph{weakly systolic groups} (i.e. groups acting geometrically by simplicial automorphisms on weakly systolic complexes) contains essentially classes of \emph{systolic groups} (i.e. groups acting geometrically by simplicial automorphisms on systolic complexes) and "cubical $CAT(-1)$ groups". For other classes of weakly systolic groups see \cite{O-sdn}.

Weakly systolic complexes are characterized as follows in \cite{O-sdn}.

\begin{de}[$SD_2^{\ast}$ property]
\label{sd2*}
For $k\geq 6$, a flag simplicial complex $X$ satisfies the $SD_2^{\ast}(k)$ property if the following two conditions hold.
\medskip

(a) $X$ does not contain $4$-wheels (cf. Section \ref{simpl}),

(b) for every $5\leq l <k$ and every $l$-wheel with a pendant triangle $\widehat W$ in $X$, there exists a vertex $v$ with $\widehat W\subseteq B_1(v,X)$.
\medskip

By \emph{$SD_2^{\ast}$ property} we mean the property $SD_2^{\ast}(6)$.
\end{de}

\begin{tw}
\label{logl}
Let $X$ be a simplicial complex satisfying the condition $SD_2^{\ast}$. Then its universal
cover is weakly systolic and in particular contractible.
\end{tw}

Weakly systolic complexes and groups possess many properties analogous to the ones of non-positively curved complexes and groups; cf. \cites{O-sdn,OCh}.

The following class of simplicial complexes will play a role in Section \ref{newa}.
A flag simplicial complex $X$ is called a \emph{complex with $SD_2^{\ast}$ links} (resp. a \emph{complex with $SD_2^{\ast}(k)$ links}) if $X$ and every of its links satisfy the property $SD_2^{\ast}$ (resp. $SD_2^{\ast}(k)$).
Observe that, by Theorem \ref{logl}, the universal cover of a complex with $SD_2^{\ast}$ links, is weakly systolic. Moreover every $6$-large simplicial complex (resp. systolic complex) is a complex with $SD_2^{\ast}$ links (resp. weakly systolic complex with $SD_2^{\ast}$ links). In fact, weakly systolic complexes with $SD_2^{\ast}$ links are weakly systolic complexes which (asymptotically) resemble the most systolic complexes.

\section{Thickening of a CAT(-1) cubical complex}
\label{prel-sd}

In this section we show how to find, for a given CAT(-1) cubical complex, an associated $5$-large, i.e. flag no-square simplicial complex.

\begin{de}[Thickening]
\label{thick}
Let $Y$ be a cubical complex. The \emph{thickening} $Th(Y)$ of $Y$ is
a simplicial complex defined in the following way. Vertices of
$Th(Y)$ are vertices of $Y$. Vertices $v_1,...,v_k$
of $Th(Y)$ span a simplex iff vertices $v_1,...,v_k\in Y$ (as
vertices of $Y$) are contained in a common cube of $Y$.
\end{de}

The following lemma is proved in \cite{O-sdn}.

\begin{lem}[Loc. $k$-large thickening]
\label{lockl}
Let $k\geq 4$ and let $Y$ be a locally $k$-large cubical complex. Then its thickening $Th(Y)$ is a locally $k$-large simplicial complex.
\end{lem}

The proof in the case $k\leq 6$ in the lemma above is much easier than the one in the general case.
Moreover the case $k=5$ is the most interesting for our construction.
Thus, for completeness and for the reader's convenience we present it below.

\begin{proof}\emph{(The case $k=4,5,6$.)}
We have to study links of vertices in $Th(Y)$. Let $v$ be a vertex.
Let, for a vertex $w\in Th(Y)_v$, the set $A_w\subseteq Y_v^{(0)}$ (here we identify the $0$-skeleton of the link of a vertex in a cubical complex with the set of vertices joined with the vertex) be the set of all vertices of $Y_v^{(0)}\subseteq Y$ belonging to the minimal cube containing $v$ and $w$.

First we prove that $Th(Y)_v$ is flag. Let $A\subseteq Th(Y)^{(0)}$ be a finite set of pairwise connected (by edges in $Th(Y)$) vertices.
Then, by the definition of the thickening we have the following.
For any two $w,w'\in A$, and for every $z\in A_w$ and $z'\in A_{w'}$, vertices $z$ and $z'$ are contained in a common cube of $Y$ (the one containing $v,w$ and $w'$) and thus $\langle z,z'\rangle \in Y_v$.
Hence $\ov A=\bigcup_{w\in A} A_w$ is a set of pairwise connected vertices in $Y_v$ and thus, by flagness of $Y_v$, it spans a simplex. It follows that $\ov A \cup \lk v \rk$ is contained in a cube of $Y$ so that $A$ is contained in the same cube and thus $A$ spans a simplex in $Th(Y)_v$. It proves that links in $Th(Y)$ are flag.

Now we prove that $Th(Y)_v$ is $k$-large, for $k=5,6$.
We do it by a contradiction. Assume there is an $l$-cycle $c=(w_0,w_1,...,w_{l-1},w_0)$ in $Th(Y)_v$ without a diagonal, for $l<k$. Then we will show that there exists an $l$-cycle $c'=(z_0,z_1,...,z_{l-1},z_0)$ in $Y_v$ without a diagonal. This would contradict $k$-largeness of $Y_v$.

Since there is no diagonal in $c$ we have that for $i\neq j$ and $i\neq j\pm 1$ (mod $l$), there exist vertices $z^i_j\in A_{w_i}$ and $z^j_i\in A_{w_j}$ not contained in a common cube (containing $v$) and thus not connected by an edge in $Y_v$. Now we treat separately two cases.
\medskip

\noindent {(\bf Case $l=4$.})
Then the cycle in $Y_v$ that leads to a contradiction is $c'=(z^0_2,z^1_3,z^2_0,z^3_1,z^0_2)$.
\medskip

\noindent {(\bf Case $l=5$.})
If $\langle z^0_2,z^3_1\rangle \in Y_v$ or $\langle z^1_3,z^4_2\rangle \in Y_v$ then, for $c'$ we can take, respectively, $(z^0_2,z^1_3,z^2_0,z^3_1)$ or $(z^1_3,z^2_4,z^3_1,z^4_2)$.
Assume that it is not true. Then if $\langle z^2_0,z^4_2\rangle \in Y_v$,
we can set $c'=(z^4_2,z^0_2,z^1_3,z^2_0,z^4_2)$, and if not we can take $c'=
(z^0_2,z^1_3,z^2_0,z^3_1,z^4_2,z^0_2)$. In every case we get a cycle $c'$ without diagonal, of length $4$ or $5$.
\end{proof}

The following lemma, whose proof is very similar to the one above, will be used in Section \ref{newa}.

\begin{lem}[$SD_2^{\ast}$ links in thickenings]
\label{sd2th}
Let $Y$ be a cubical complex with links satisfying the $SD_2^{\ast}$ property. Then links in $Th(Y)$ satisfy the $SD_2^{\ast}$ property as well.
\end{lem}
\begin{proof}
By Lemma \ref{lockl} it follows that $Th(Y)$ is locally $5$-large, thus satisfying the condition (a) from the definition of the $SD_2^{\ast}$ property (Definition \ref{sd2*}).

Hence it leaves to establish the condition (b). Let $\widehat W=(v_0;v_1,...,v_5;t)$ be a $5$-wheel with a pendant triangle in $Th(Y)_v$, for some vertex $v\in Y$.

Then, since the cycle $(v_1,v_2,...,v_5,v_1)$ has no diagonal, we can find a cycle $(w_1,w_2,...,w_5,w_1)$ without a diagonal, with $w_i\subseteq A_{v_i}$ (here we use the notation from the proof of Lemma \ref{lockl})---this can be done in the same way as in the proof of Lemma \ref{lockl}.
Furthermore, we can extend this cycle to a $5$-wheel with a pendant triangle $\widehat W'=(w_0;w_1,...,w_5;s)$ in $Y_v$, such that $w_i\in A_{v_i}$ and $s\in A_{t}$. Then, by our assumptions on links in $Y$, there exists a vertex $z\in Y_v$ with $\widehat W'\subseteq B_1(z,Y_v)$. It follows that $\widehat W \subseteq B_1(z,Th(Y)_v)$ and hence the lemma follows.
\end{proof}

Now, using the local results above, we prove that in the case of simply connected locally $5$-large cubical complexes (i.e. CAT(-1) cubical complexes) their thickenings are $5$-large. The (versions of the) following two lemmas were proved also in \cite{O-sdn}.

\begin{lem}[Thickening of CAT(0) c.c.]
 \label{flag thick}
Let $Y$ be a simply connected locally $4$-large cubical complex. Then
$Th(Y)$ is a $4$-large (i.e. flag) simplicial complex.
\end{lem}
\begin{proof}
Let $v_0,v_1,...,v_k\in Th(Y)$ be vertices pairwise joined by edges
in $Th(Y)$. We have to prove that they span a simplex in $Th(Y)$, that
means they are contained (as vertices of $Y$) in a common cube.

Let $B_i=B_1(v_i,Th(Y))$. Since the balls $B_i$ pairwise intersect, by \emph{Helly property} \cite[Proposition 2.6]{BaVe} they all intersect. Let $v\in \bigcap_i B_i$.
Since $Y_v$ is flag, by Lemma \ref{lockl} we have that $Th(Y)_v$
is flag. Since $v_0,v_1,...,v_k\in B_1(v,Th(Y))$ we get that these
vertices span a simplex in $Th(Y)$.
\end{proof}

\begin{prop}[Thickening of $CAT(-1)$ c.c.]
\label{5l thick}
Let $Y$ be a simply connected locally $5$-large cubical complex. Then
$Th(Y)$ is a $5$-large simplicial complex.
\end{prop}
\begin{proof}
By Lemma \ref{flag thick} we have that $Th(Y)$ is flag. Thus it only leaves to show that every $4$-cycle in $Th(Y)$ has a diagonal. Let $(v,w,z,w',v)$ be a cycle in the $1$-skeleton of $Th(Y)$.
Let $c$ and $c'$ be maximal cubes containing, $z$ and, respectively,
$w$ and $w'$. Let $B$ be a subgraph of $Y^{(1)}$ spanned by
$B_1(v,Th(Y))^{(0)}$. Then it follows from \cite[Proof of Proposition 2.6]{BaVe} that $c^{(1)}$, $c'^{(1)}$ and $B$ are convex subgraphs of $Y^{(1)}$. Thus by \cite[Lemma 4.2]{Ch4}, they are the so called \emph{gated} subgraphs. By \cite[Section 1]{BaCh} gated subgraphs satisfy the Helly property thus there exists a vertex $u$ in their intersection.
Hence we have that $v,w,w',z \in Th(Y)_u$ and, by Lemma \ref{lockl}, we have that the cycle $(v,w,z,w',v)$ has a diagonal since $Th(Y)$ is locally $5$-large.
\end{proof}

Actually, in \cite{O-sdn} we prove the following stronger version of the above result. It follows quite easily from Theorem \ref{logl} and will be not used in the basic construction.

\begin{prop}[Weakly systolic thickening of $CAT(-1)$ c.c.]
\label{5l wsys}
Let $Y$ be a simply connected locally $5$-large cubical complex. Then
$Th(Y)$ is a weakly systolic $SD_2^{\ast}(k)$ complex, for every $k\geq 6$.
\end{prop}

\begin{lem}
\label{bjorner}
Let $Y$ be a finite dimensional cubical complex. Then $Th(Y)$ is homotopically equivalent with $Y$.
\end{lem}
\begin{proof}
It follows immediately from the Nerve Theorem of Borsuk \cite{Bor}.
\end{proof}

\begin{lem}
\label{systh}
Let $Y$ be a locally $k$-large, simply connected cubical complex, for some $k\geq 4$. Then $Th(Y)$ is a $k$-systolic, contractible simplicial complex.
\end{lem}
\begin{proof}
By Lemma \ref{lockl}, $Th(Y)$ is locally $k$-large, and by Lemma \ref{bjorner}, it is contractible (since CAT(0) cubical complexes are contractible).
\end{proof}

For a cubical complex $Y$ we can identify links of vertices in $Y$ with subcomplexes in links of vertices of $Th(Y)$. The following lemma is obvious. It will be used in Section \ref{new}.

\begin{lem}
\label{fullsbc}
Let $v$ be a vertex of a locally flag cubical complex $Y$. Then every full subcomplex of $Y_v$ is a full subcomplex of $Th(Y)_v$.
\end{lem}

\rem
In a more general setting of cell complexes, the thickening was invented by T. Januszkiewicz and, independently, by the author; cf. \cite{O-sdn}.
In the case of cubical complexes, a construction similar to the thickening has been used in graph theory. For a graph $G$ being the $1$-skeleton of a cubical complex $Y$ (such graphs are called \emph{median graphs}), a graph $G^{\Delta}$ is defined as the $1$-skeleton of $Th(Y)$; cf. \cite{BaCh} (we use their notation here).

\section{The basic construction}
\label{proof}
In this section we present in details the construction of high dimensional Gromov hyperbolic right angled Coxeter groups as described roughly in the Introduction---Section \ref{intro}. We will keep here the notation introduced there.

By the Moussong's theorem (Theorem \ref{mous}), a Gromov hyperbolic right angled Coxeter group $W$ is determined by its nerve $L=L(W,S)$ (where $S$ is a given set of generators of $W$) that is itself a finite $5$-large simplicial complex. By Theorem \ref{vcdcox}, the virtual cohomological dimension of $W$ is estimated by a "dimension" of $L$.
Thus the question reduces to a construction of finite $5$-large simplicial complexes $X$ of arbitrarily large "dimension", i.e. with $H^n(X)\neq 0$ for large $n$.

We construct such complexes by induction. Their "dimension" will increase at each step.

The base of the induction is a finite $5$-large simplicial complex $X_0$ with $H^{n_0}(X_0)\neq 0$. As an example one can choose $X_0$ to be a finite graph of girth at least $5$ (i.e. not containing cycles of length less than $5$) containing a $k$-cycle for some $k\geq 5$.

The induction step is as follows. The input data is a finite $5$-large simplicial complex $X$ with
$H^{n}(X)\neq 0$, for $n\geq n_0$.
As an output we obtain a finite $5$-large simplicial complex $X'$ with $H^{n+1}(X')\neq 0$.

Now we present the induction step in details following the overview given in Introduction. It consists itself of the three following steps.
\medskip

\noindent
{\bf Step 1.}
Let $(W,S)$ be the right-angled Coxeter system whose nerve is $X$, i.e. $X=L(W,S)$. By Theorem \ref{vcdcox}, the virtual cohomological dimension of $W$ is $\mr{vcd}\; W\geq n+1$. The Davis complex $\Sigma =\Sigma(W,S)$ is a contractible cubical complex in which links of all vertices are isomorphic with $X$. Thus it is a contractible locally $5$-large
cubical complex.
\medskip

\noindent
{\bf Step 2.}
In this step we find a torsion-free subgroup $H<W$ of large finite index  such that $Y=H\backslash \Sigma$ is a locally $5$-large cubical complex with $H^{n+1}(Y)\neq 0$. The only difficulty here is to prove the last inequality. We do it, using standard tools (cf. \cite{Dav-b}), by constructing a nontrivial cohomology class in $H^{n+1}(Y)$; see Lemma \ref{f'<>0} below. We use simplicial cohomology w.r.t. the standard triangulation of $\Sigma$, in which the Davis chamber $K$ is the cone over the barycentric subdivision of the nerve $L$.
\medskip

For $H< W$, let $Y=H\backslash \Sigma$ and let $p\colon \Sigma \to Y$ be the quotient map.
Let $\mathcal K$ be the \emph{set of copies of $ K$ in $Y$} i.e.
${\mathcal K} =\lk  p(w K)|w\in  W \rk$.

\begin{de}[Orientation]
\label{orient}
An \emph{orientation} of $Y=H\backslash \Sigma$ is a map
$\epsilon \colon \mathcal K \to \lk -1,1 \rk$ such that $\epsilon(p({s} w
K))=-\epsilon(p( w K))$, for every $ w\in  W$ and $ s\in  S$.
\end{de}

Observe that $\epsilon ( w  K)=(-1)^{l_{ S}( w)}$ (here $l_{ S}( w)$ denotes the length of $ w$ in the word metric on $ W$ w.r.t. $ S$) defines an orientation on $Y=\Sigma$, i.e. in the case $H=\lk 1 \rk$.

By Corollary D.1.4 of \cite{Dav-b} Coxeter groups are virtually torsion-free and thus there exists a finite index torsion-free subgroup $H$ of
$ W$. Coxeter groups are residually finite (cf. e.g. \cite[Section 14.1]{Dav-b}) so that we
can
choose $H$ whose minimal displacement is "big". For our purposes we make the following choice.
We take a torsion-free subgroup $H< W$ such that for the action (induced by the $W$-action on $\Sigma$) of $H$ on $Th(\Sigma)$ (on the thickening of $\Sigma$---cf. Section \ref{prel-sd}) we have $\inf \lk \mr d(v,hv)\; | \; v\in (Th(\Sigma))^{(0)},\; h\in H \rk \geq 5$, where $\mr d(\cdot , \cdot)$ is the usual metric on $(Th(\Sigma))^{(1)}$.

As in the case of manifolds which have two-sheeted orientable covering, the same is true here, with an analogous proof.

\begin{lem}
\label{orient cov}
There exists a two-sheeted covering of $Y$, which admits an orientation.
\end{lem}

Thus w.l.o.g. we can assume that there is an orientation $\epsilon \colon {\mathcal K}\to \lk-1 ,1\rk$ of $Y$.
This corresponds to replacing $H$ by its subgroup of index two.
It is clear that then there exists an orientation of $\Sigma$,
$\epsilon'\colon \lk  w K|  w\in  W \rk \to \lk -1, 1 \rk$ such that
$\epsilon(p( w K))=\epsilon'( w( K))$. We will write $\epsilon$ instead of $\epsilon'$---it
should not lead to misunderstanding.

\begin{de}[Antisymmetrization of cochains]
\label{antisym}
Let $\epsilon$ be an orientation of $Y$. For $h\in C^{\ast}(Y)$ let
$$
a(\epsilon,h,\sigma)=\sum _{ K'\in {\mathcal K}}\epsilon(
K')h(p({ w}_{ K'}\sigma)),
$$
where $\sigma$ is a simplex in $ K$, ${ w}_{ K'}\in  W$ and
$p({ w}_{ K'}\sigma)\in  K'$.
The \emph{antisymmetrization wrt an orientation $\epsilon$} of a cochain $h$ is a cochain  $a_{\epsilon}(h)$ defined as follows:
$$
a_{\epsilon}(h)(p({ w}\sigma))=\frac{\epsilon({ w} K)}{|{\mathcal K}|}a(\epsilon,h,\sigma).
$$
where ${ w}\in  W$ and $\sigma$ is a simplex in $ K$.
\end{de}

\begin{lem}[Properties of $a_{\epsilon}(h)$]
\label{prop a}
Let $h\in C^{\ast}(Y)$. The antisymmetrization has the following properties:
\begin{enumerate}
\item $\delta a_{\epsilon}(h)=a_{\epsilon}(\delta h)$,

\item $a_{\epsilon}(a_{\epsilon}(h))=a_{\epsilon}(h)$.
\end{enumerate}
\end{lem}
\begin{proof}
Let ${ w}\in  W$ and let $\sigma$ be a simplex in $ K$. Then we have the following.

\noindent
(1)
\begin{align*}
\delta a_{\epsilon}(h)(p({ w}\sigma)) & = a_{\epsilon}(h)(\partial p({ w}\sigma))=a_{\epsilon}(h)(p({ w}\partial \sigma))=\frac{\epsilon({ w} K)}{|{\mathcal K}|}a(\epsilon,h,\partial \sigma) \\
& = \frac{\epsilon({ w} K)}{|{\mathcal K}|}
\sum _{ K'\in {\mathcal K}}\epsilon(
K')h(\partial (p ({ w}_{ K'}\sigma))) \\
& =\frac{\epsilon({ w} K)}{|{\mathcal K}|}
\sum _{ K'\in {\mathcal K}}\epsilon(
K')\delta h(p({ w}_{ K'}\sigma)) \\
& =\frac{\epsilon({ w} K)}{|{\mathcal K}|}a(\epsilon,\delta h,\sigma)=
a_{\epsilon}(\delta h)(p({ w}\sigma)).
\end{align*}

\noindent
(2)
\begin{align*}
a_{\epsilon}(a_{\epsilon}(h))(p({ w}\sigma)) & = \frac{\epsilon({ w} K)}{|{\mathcal K}|}a(\epsilon,a_{\epsilon}(h),\sigma)=
\frac{\epsilon({ w} K)}{|{\mathcal K}|}\sum _{ K'\in {\mathcal K}}\epsilon( K')a_{\epsilon}(h)(p({ w}_{ K'}\sigma)) \\
& = \frac{\epsilon({ w} K)}{|{\mathcal K}|}
a(\epsilon,h,\sigma)
\sum _{ K'\in {\mathcal K}}\epsilon( K')\frac{\epsilon( K')}{|{\mathcal K}|}=a_{\epsilon}(h)(p({ w}\sigma)).
\end{align*}

\end{proof}

By our inductive assumptions, there exists a cocycle $f\in Z^{n+1}(K,K^{S})$ representing a non-trivial class
$[f]\in H^{n+1}( K,K^{S})=H^n(X)$. We can treat $f$ as a cochain on $Y$ (with support in $p(K)$) and define a cochain $ f'=|\mathcal K|a_{\epsilon}( f)$ on $Y$.

\begin{lem}[Non-triviality of ${[ f']}$]
\label{f'<>0}
The cochain $ f'$ is a cocycle representing a non-trivial class $[ f']\in H^{n+1}(Y)$.
\end{lem}
\begin{proof}
First we prove that $\delta  f'=0$.
this follows from Lemma \ref{prop a}, since $\delta  f'=\delta |\mathcal K|a_{\epsilon}( f)=|\mathcal K|a_{\epsilon}(\delta  f)=0$ (the last equality follows from the fact that $f$ vanishes on $K^S$).

Now we prove that $[ f']\neq 0$ in $H^{n+1}(Y)$. Assume this is not so.
Let $ f'=\delta g'$ for some cochain $g'\in C^{n}(Y)$.
Then, by Lemma \ref{prop a}, $\delta a_{\epsilon}(g')=a_{\epsilon}(\delta g')=a_{\epsilon}( f')=a_{\epsilon}(|\mathcal K|a_{\epsilon}( f))= f'$.
Observe that $a_{\epsilon}(g')$ vanishes on $p(wK^{ S})$, for every $w\in W$.
Let $g\in C^{n}(K,K^S)$ be defined as $g(\sigma)=\frac{1}{|\mathcal K|}a(\epsilon, a_{\epsilon}(g'), \sigma)$, for $\sigma$ in $K$.
Then, for any $(n+1)$-simplex $\sigma$ in $K$, we have
\begin{align*}
\delta g(\sigma) & =\frac{1}{|\mathcal K|}a(\epsilon, a_{\epsilon}(g'), \partial \sigma) =
\frac{1}{|\mathcal K|}\sum _{ K'\in {\mathcal K}}\epsilon(
K')a_{\epsilon}(g')( p ({ w}_{ K'}\partial \sigma)) \\
& =\frac{1}{|\mathcal K|}
\sum _{ K'\in {\mathcal K}}\epsilon(
K')\frac{\epsilon(K')}{|{\mathcal K}|}a(\epsilon, g', \partial \sigma)
=\frac{1}{|\mathcal K|}\sum _{ K'\in {\mathcal K}}{\epsilon(K')}g'(p ({ w}_{ K'}\partial \sigma))\\
& =\frac{1}{|\mathcal K|}\sum _{ K'\in {\mathcal K}}{\epsilon(K')}f'(p ({ w}_{ K'} \sigma))=
\frac{1}{|\mathcal K|}|\mathcal K|
\sum _{ K'\in {\mathcal K}}{\epsilon(K')}a_{\epsilon}(f)(p ({ w}_{ K'} \sigma)) \\
& =
\sum _{ K'\in {\mathcal K}}{\epsilon(K')}\frac{\epsilon(K')}{|{\mathcal K}|}a(\epsilon, f, \sigma)=\sum _{ K'\in {\mathcal K}}{\epsilon(K')}
f(p ({ w}_{ K'} \sigma))=f(\sigma).
\end{align*}
This contradicts the fact that $[f]\neq 0$ in $H^{n+1}(K,K^S)$ and finishes the proof of the lemma.
\end{proof}
\medskip

\noindent
{\bf Step 3.}
Let $X'=Th(Y)$. Then, by Lemma \ref{bjorner} and Lemma \ref{f'<>0}, $H^{n+1}(X')=H^{n+1}(Y)\neq 0$. Observe that $X'=H\backslash Th(\Sigma)$. Since, by Lemma \ref{5l thick}, the complex $Th(\Sigma)$ is $5$-large, our choice of $H$ (with large minimal displacement) guarantees that $X'$ is $5$-large.

Let us conclude what we get in Steps 1--3, in the following theorem (cf. Main Theorem from Introduction).

\begin{tw}
\label{maint}
Let $X$ be a finite $5$-large simplicial complex such that $H^n(X)\neq 0$.
Then the complex $X'$ (as obtained in Steps 1--3) is a finite $5$-large simplicial complex and $H^{n+1}(X')\neq 0$. Thus the right-angled Coxeter group with nerve $X'$ is Gromov hyperbolic of virtual cohomological dimension at least $n+2$.
\end{tw}

\section{New examples of hyperbolic Coxeter groups}
\label{new}
In this section we show that our construction provides also examples of highly dimensional Gromov hyperbolic Coxeter groups that are very (i.e. asymptotically) different from the ones constructed before. In fact, all such Coxeter groups obtained by using constructions from \cites{JS1,JS2,H,AB+} are systolic.
Systolic groups are in a way two-dimensional---they "do not contain asymptotically" spheres above dimension one; cf. \cites{JS3, O-ciscg, O-ib7scg,OS}.
We show how to get non-systolic examples (cf. Corollary \ref{nonsys} below), and moreover, we point out that groups obtained by our construction can "contain asymtotically" spheres up to dimension $3$ (cf. Remarks after Corollary \ref{nonsys}).

\begin{prop}[Subgroups]
\label{sbgp}
Let $Z$ be a full subcomplex of a finite $5$-large simplicial complex $X_0$. Let $X$ be a finite $5$-large simplicial complex obtained by the basic construction starting with $X_0$ (i.e. by performing several times Steps 1--3 from Section \ref{proof}) and let $(W,S)$ be the resulting Coxeter system (i.e. $X=L(W,S)$). Then the right-angled Coxeter group $W'$ with the nerve $Z$ (i.e. $Z=L(W',S')$) is a subgroup of $W$.
\end{prop}
\begin{proof}
By repeated use of Lemma \ref{fullsbc}, we get that $Z$ is a full subcomplex of $X$. Thus it corresponds to the subgroup $W'$ of $W$; cf.   \cite[Section 8.8]{Dav-b}.
\end{proof}

\begin{cor}[Non-systolic]
\label{nonsys}
   If $X_0$ is the nerve of a group $W'$ which is not asymptotically hereditarily aspherical (AHA; cf. \cites{JS3,OS}), then a group $W$ obtained by the basic construction is not systolic (or not acting geometrically on a weakly systolic complex with $SD_2^{\ast}$ links).
\end{cor}
\begin{proof}
It follows directly from Proposition \ref{sbgp} above and from the fact that subgroups of systolic groups (or of groups acting geometrically on a weakly systolic complexes with $SD_2^{\ast}$ links) are AHA; cf. \cites{JS3,OS}.
\end{proof}

\noindent
{\bf Examples.} Here we show how our construction can lead to the non-systolic examples. Let $X_0$ be a nerve of a right-angled Coxeter group $W'$ acting geometrically on the $k$-dimensional real hyperbolic space $\mathbb{H}^k$ (such groups exists only for $k=1,2,3,4$; cf. \cite{V}). Then $W'$ is not AHA (cf. \cites{O-ciscg,OS}) and thus $W$ (a group obtained by our construction) is not AHA, in particular not systolic (and do not act geometrically on weakly systolic complex with $SD_2^{\ast}$ links); cf. \cites{JS3,O-sdn,OS}.

\medskip

\rems
1) A group $W$ as in the example above "contains
$\partial W'=\mathbb{S}^{k-1}$ at infinity". AHA groups, and in particular systolic groups (or groups acting geometrically on weakly systolic complexes with $SD_2^{\ast}$ links), do not "contain asymptotically" spheres above dimension $1$; cf. \cites{JS3,O-ciscg,OS}. 
For other relevant asymptotic invariants that can distinguish our groups from e.g. systolic groups see \cites{O-ib7scg,OS}.

2) We do not know whether Gromov hyperbolic Coxeter group can "contain asymptotically" spheres above dimension $3$. A (form of a) conjecture by Januszkiewicz-\' Swi\c atkowski says that the Gromov boundary of a simply connected locally $5$-large cubical complex (i.e. CAT(-1) cubical complex) cannot contain a sphere of dimension above $3$. Since hyperbolic Coxeter groups act geometrically on such cubical complexes (their Davis complexes), this conjecture implies that their boundaries cannot contain highly dimensional spheres.

3) Yet another way of distinguishing our groups from systolic groups 
(or groups acting geometrically on a weakly systolic complexes with $SD_2^{\ast}$ links) can rely on the following fact. A finitely presented subgroup of a torsion-free $k$-systolic group (resp. a group acting geometrically on a weakly systolic complex with $SD_2^{\ast}(k)$ links)
is $k$-systolic (resp. acts geometrically on a weakly systolic complex with $SD_2^{\ast}(k)$ links); cf. \cite{Wis} for the systolic case and \cite{O-sdn} for the weakly systolic case.

\section{New examples of groups with some asphericity properties}
\label{newa}

This section complements in a way the previous one. We show namely that our construction leads to (a priori new) systolic groups and groups having interesting asymptotic properties exhibited by systolic groups as in e.g. \cites{JS3, O-ciscg, O-ib7scg}.

\begin{prop}[Systolic groups]
\label{syst}
Let $k\geq 6$ and let $X$ be a finite $5$-large simplicial complex obtained by the basic construction (as in Section \ref{proof}). Assume that the two following sets of conditions are satisfied.
\medskip

\noindent
1) The initial complex $X_0$ is $k$-large (cf. Section \ref{proof});
\medskip

\noindent
2) At every step of the basic construction the group $H$ (cf. Step 2 in Section \ref{proof}) is chosen so that the minimal displacement of the $H$-action on $Th(\Sigma)^{(1)}$ is at least $k$.
\medskip

\noindent
Then the resulting complex $X$ is $k$-large and, consequently, the Coxeter group $W$ with the nerve $X$ is $k$-systolic.
\end{prop}
\begin{proof}
If the complex $X$ (in the inductive step of the basic construction) is $k$-large then $\Sigma$ is locally $k$-large and, by Lemma \ref{systh}, its thickening $Th(\Sigma)$ is $k$-systolic. By the assumption 2, we have that then $X'=H\backslash Th(\Sigma)$ is $k$-large. Thus the lemma follows inductively.
\end{proof}

\begin{prop}[Complexes with $SD_2^{\ast}$ links]
\label{sd*}
Let $X$ be a finite $5$-large simplicial complex obtained by the basic construction with the initial complex $X_0$ being a complex with $SD_2^{\ast}$ links. Then the resulting complex $X$ is a complex with $SD_2^{\ast}$ links and, consequently, the Coxeter group $W$ with the nerve $X$ acts on a weakly systolic complex with $SD_2^{\ast}$ links.
\end{prop}
\begin{proof}
If the complex $X$ (in the inductive step of the basic construction) is a complex with $SD_2^{\ast}$ links then links in $\Sigma$
satisfy trivially the property $SD_2^{\ast}$. By Lemma \ref{sd2th}, it follows that links in $Th(\Sigma)$ and thus also in $X'=H\backslash Th(\Sigma)$ satisfy the property $SD_2^{\ast}$.
Finally, by Lemma \ref{5l wsys}, $Th(\Sigma)$ and thus also $X'$ satisfy themselves the property $SD_2^{\ast}$. Hence the proposition follows inductively.
\end{proof}

\begin{cor}[Asphericity properties]
\label{asph}
A group $W$ as in Proposition \ref{sd*} has the following properties.
\begin{enumerate}
\item $W$ is not simply connected at infinity.

\item Higher (then the first) homotopy pro-groups at infinity vanish for $W$.

\item The Gromov boundary of $W$ is strongly hereditarily aspherical (i.e. SHA; cf. e.g. \cite{O-ib7scg}).

\end{enumerate}
\end{cor}
\begin{proof}
(1) and (2) follow from the fact that the same properties hold for all groups acting geometrically on weakly systolic complexes with $SD_2^{\ast}$ links; cf. \cite{O-sdn}.

(3) follows from the fact that the Gromov boundary of a simply connected locally $5$-large cubical complex with $SD_2^{\ast}$ links is SHA; cf. \cite{O-sdn}.
\end{proof}

\rems 1) Corollary \ref{asph} yields that our construction provides (a priori) new examples of, e.g. hyperbolic groups with SHA boundary. Before, the only known examples of such groups in higher dimensions were some systolic groups; cf. \cites{O-ib7scg, O-sdn}. Other examples can be obtained by taking finitely presented subgoups---finitely presented subgroups of fundamental groups of finite complexes with $SD_2^{\ast}$ links are of the same type; cf. \cite{O-sdn}

2) In \cite{OS} we prove that groups acting geometrically on weakly systolic complexes with $SD_2^{\ast}$ links are asymptotically hereditarily aspherical (shortly AHA); cf. \cite{JS3}. Again, Proposition \ref{sd*}, yields (a priori) new examples of highly dimensional
AHA groups. The only other ones are systolic groups.
Actually, it is an interesting phenomenon that there exist AHA groups or groups with SHA boundary in dimensions above two.

3) Asymptotic properties listed in Corollary \ref{asph} hold also for systolic groups (cf. Section \ref{weak}).
It should be noticed here that these properties does not hold for many classical groups. In particular they distinquish our groups from, e.g.: fundamental groups of closed manifolds covered by $\mathbb R^n$ (in particular non-positively curved manifolds), with $n\geq 3$; Coxeter groups acting on $\mathbb R^n$, $n\geq 3$; and some lattices in isometry groups of related buildings; cf. \cites{JS3, O-ciscg, OS}.
For other relevant properties of groups acting geometrically on weakly systolic complexes with $SD_2^{\ast}$ links see \cites{O-sdn,OS}.

\section{Final remarks}
\label{final}

\subsection{Variants of the basic construction}
\label{var}
There are many variants of the basic construction corresponding to the choice of the complex $Z$ playing the role of $\Sigma$ in Step 1.(cf. Section \ref{proof}), for a fixed $X$.
Let us describe here one of them.

For a given $X$, $W$ and $\Sigma$ as in Step 1. of the basic construction, let $Z$ be a right angled building associated with $W$. A simple construction of such buildings is described in \cite{DO}. It is also showed there (reproving known facts in a new way) that groups of isometry of such buildings admit uniform lattices. In \cite{O-sdn} it is shown that for any such building there exists a locally $5$-large simplicial complex $Th(Z)$, being an analogue of the thickening of a cubical complex described in Section \ref{prel-sd}. Thus, for an appropriate $H$ being a lattice in $\mr{Isom}(Z)$, the quotient $H\backslash Th(Z)$ is a finite $5$-large simplicial complex of cohomological dimension at least $n+1$ (since $Z$ contains copies of $\Sigma$). And thus we can continue the induction.

\subsection{Non right-angled case}
\label{nra}

Having a Gromov hyperbolic right-angled Coxeter group $W$ with a presentation $W=\langle S \;|\; (st)^{m_{st}}\rangle$, one can construct a non right-angled Coxeter group $W^f=\langle S \;|\; (st)^{m^f_{st}}\rangle$ as follows.
If $m_{st}\neq \infty$, then $m^f_{st}=m_{st}$, otherwise $m^f_{st}$ is an arbitrary number larger than $4$, or $\infty$. It is proved in \cite{JS1} that $W^f$ is then Gromov hyperbolic and that $\mr {vcd}(W^f)=\mr {vcd}(W)$, if $\mr {vcd}(W)\geq 2$.

\subsection{Parametrization of the construction}
\label{par}

Combinatorial balls in $Th(\Sigma)$ (for $\Sigma$ being the Davis complex of a right-angled Coxeter group, as in Step 1. in Section \ref{proof}) are fundamental domains for finite index subgroups of $W$ (as in Step 1. of the basic construction). Thus sequences $\ov s$ of the following form correspond to Coxeter groups obtained by the basic construction and having dimension at least $(n_0+k)$. Let $X_0$ be a finite $5$-large simplicial complex of cohomological dimension at least $n_0$. Then let $\ov s=(X_0;r_1,r_2,...,r_k)$, where $r_i\geq 5$ (a radius of a combinatorial ball in the thickening of the Davis complex at the $i$-th step of the construction).

Moreover, if one wants to get $k$-systolic groups then it is enough to assure that $X_0$ is $k$-large and that $r_i\geq k$; cf. Proposition \ref{syst}.



\end{document}